\documentclass[11pt,a4paper,reqno]{amsart}

\usepackage{amsmath}
\usepackage{amsfonts}
\usepackage{amssymb}
\usepackage{latexsym}
\usepackage{graphics}

\numberwithin{equation}{section}
\newtheorem{proposition}{Proposition}[section]
\newtheorem{lemma}[proposition]{Lemma}

\newtheorem{theorem}[proposition]{Theorem}
\newtheorem{conjecture}[proposition]{Conjecture}
\theoremstyle{definition}

\newtheorem{example}[proposition]{Example}

\newcommand{\sn}{S_n}
\newcommand{\F}{\mathcal{F}}
\newcommand{\T}{\mathcal{T}}
\newcommand{\A}{\mathcal{A}}
\newcommand{\R}{\mathcal{R}}

\renewcommand{\H}{\mathcal{H}}
\newcommand{\X}{\mathcal{X}}
\newcommand{\Y}{\mathcal{Y}}
\newcommand{\Pair}{\mbox{\scshape{Pair}}}

\newcommand{\Refine}{\mbox{\scshape{Refine}}}
\newcommand{\Expand}{\mbox{\scshape{Expand}}}
\newcommand{\Random}{\mbox{\scshape{Filter}$^{\star}$}}
\newcommand{\Free}{\mbox{\scshape{Free}}}

\newcommand{\GETOUT}[1]{}

\newcommand{\preimage}[1]{N({#1})}
\newcommand{\myrest}[2]{{#1}[{#2}]}
\newcommand{\property}{\beta}

\title{Concerning the shape of a geometric lattice}
\author{W. M. B. Dukes}
\email{mark.dukes@ucd.ie}
\thanks{Supported by EC's Research Training Network `Algebraic Combinatorics in Europe', grant 
HPRN-CT-2001-00272 while the author was at Universit\`{a} di Roma Tor Vergata, Italy and Universit\'{e} Bordeaux 1, France.}
\address{School of Mathematical Sciences, University College Dublin, Ireland.}

\begin{document}
\maketitle
\begin{abstract}
A well-known conjecture states that the Whitney numbers of the second kind of a geometric lattice (simple matroid) are logarithmically concave. 
We show this conjecture to be equivalent to proving an upper bound on the number of new copoints in the free erection of the associated simple matroid $M$. 
A bound on the number of these new copoints is given in terms of the copoints and colines of $M$.
Also, the points-lines-planes conjecture is shown to be equivalent to a problem concerning the number subgraphs of a
certain bipartite graph whose vertices are the points and lines of a geometric lattice.
\end{abstract}

\section{Introduction}
In this section we introduce some standard terminology concerning matroids and their erections,
state conjectures concerning the shape of a geometric lattice and highlight results to date.
In Section 2 we define the free erection of a matroid via constructs from Knuth's~\cite{knuthrandom} random matroids algorithm,
and prove in Theorem~\ref{theorem1} that the free erection of a simple matroid is the erection containing
maximum number of new copoints.
The relevance of the free erection to the log-concavity conjecture is made clear in Theorem~\ref{lcequiv},
while in Theorem~\ref{freebound} we give an upper bound on the number of copoints in the free erection of geometric lattice.
In Section 3, using a recent result from Dukes~\cite{dukes.slc}, we give an equivalent formulation of the points-lines-planes conjecture in terms 
of the number of subgraphs of a bipartite graph.

Let $M$ be a simple matroid on the $n$-element set $S_n$ with rank function $\rho$. 
The rank-$i$ flats are denoted by $\F_i(M)$ for all $0\leq i \leq \rho(M)$.
The collections of {\it{copoints}}, {\it{colines}} and {\it{coplanes}} of $M$ are denoted $\F_{\rho(M)-1}(M)$, $\F_{\rho(M)-2}(M)$ and $\F_{\rho(M)-3}(M)$, respectively.
The flats of a simple matroid, ordered by inclusion, form a geometric lattice.
Where there is no confusion, $\F_i$ will be used instead of $\F_i(M)$. 
The numbers $W_i(M):=|\F_i(M)|$ are called the {\it{Whitney numbers (of the 2$^{\mbox{\tiny nd}}$ kind)}} of the geometric lattice $M$.
The $k$-{\it{truncation}} of $M$ is the rank-$(k+1)$ matroid $U_k(M)$ with flats $\F(U_k(M)) := \{S_n\} \cup \{F \in \F_i(M): 0\leq i \leq k\}$.
The {\it{closure}} of a set $A\subset S_n$ is the intersection of all flats in $\F$ containing $A$.
A set $A$ is {\it{$k$-closed}} if it contains the closures of all its $j$-element subsets, for all $j\leq k$.

For matroids $M$ and $N$ on $S_n$, call $N$ an {\it{erection}} of $M$ if
the flats of $M$ are precisely those flats of $N$ not of rank $\rho(M)$.
Allow $M$ to be an erection of itself, called the {\it{trivial erection}}.
Notice that a matroid $N$ is a non-trivial erection of $M$ if and only if $M$ is the $\rho(M)-1$ truncation of $N$. 

\begin{conjecture}[Mason~\cite{mason}]
\label{lcstatement}
Let $M$ be a simple matroid with $\rho(M)>2$, then for all $0<k<\rho(M)$,
\renewcommand{\theenumi}{\roman{enumi}}
\begin{enumerate}
\item $W_k(M)^2 \geq W_{k-1}(M) W_{k+1}(M)$,\\[-0.2em]
\item $W_k(M)^2 \geq \dfrac{k+1}{k} W_{k-1}(M) W_{k+1}(M)$,
\item $W_k(M)^2 \geq \dfrac{(k+1)(n-k+1)}{k(n-k)} W_{k-1}(M) W_{k+1}(M)$.
\end{enumerate}
\renewcommand{\theenumi}{\arabic{enumi}}
\end{conjecture}
Of course (iii) $\Rightarrow$ (ii) $\Rightarrow$ (i).
The idea behind the coefficient in (iii) is that the maximum of the ratio $W_k(M)^2/W_{k-1}(M) W_{k+1}(M)$ is 
thought to be attained when $W_i(M)={n\choose i}$.
The particular case of (iii) with $k=2$ is known as the `points-lines-planes conjecture'. 
(The term `points-lines-planes' conjecture was originally coined in Welsh~\cite[p.289]{welshbook} for (ii) with $k=2$, 
however Seymour's more recent paper deals with the more general inequality (iii) with $k=2$.)
The above conjectures are strengthenings of Rota's~\cite{rota64} conjecture that the Whitney numbers of a matroid are unimodal.
It is well known that $W_1(M) \leq W_k(M)$ for all $k>1$, but recently Kung~\cite{kung} has shown that 
$W_2(M) \leq W_3(M)$ for simple matroids $M$ of rank greater than 5 and in which all lines contain the same number of points.
Stonesifer~\cite{stone} showed the points-lines-planes conjecture to be true for all graphic matroids.
Seymour~\cite{seymour} generalized the result and proved the conjecture to be true for all simple matroids $M$ such that $|F|\leq 4$, for all $F \in \F_2(M)$
(i.e. no lines contain more than four points).
His `localized' proof relies on exhibiting a function $k(x,y)$ which satisfies a collection of
inequalities relating the number of lines to the number of planes.

In what follows we deal with the log-concavity conjecture (i). 
An almost trivial fact which ties up both ends of (\ref{lcstatement}) (i) is 
\begin{lemma}
\label{endclause}
Let $M$ be a simple matroid on $S_n$ of rank $r\geq 2$. Then $W_1(M)^2 \geq W_0(M) W_2(M)$ and $W_{r-1}(M)^2 \geq W_{r-2}(M) W_r(M)$.
\end{lemma}
\begin{proof}
The first inequality follows since $W_0(M)=1$, $W_1(M) = n$ and $W_2(M) \leq {n\choose 2}$.
The second inequality holds because $W_r(M)=1$ and the collection of intersections of all copoints in $\F_{r-1}$ generates all flats
of lower rank, which includes the collection $\F_{r-2}$.
\end{proof}

Erections of matroids (geometries) were introduced by Crapo\cite{crapo70i,crapo70ii}, 
where the following fundamental result was first proved. 
\begin{theorem}[Crapo\cite{crapo70ii}]
A set $\H$ of subsets of $S_n$ is the set of copoints of an erection of $M$ if and only if
\renewcommand{\theenumi}{\roman{enumi}}
\begin{enumerate}
\item each $H\in \H$ spans $M$
\item each $H\in \H$ is $(\rho(M)-1)$-closed
\item each basis for $M$ is contained in a unique $H\in \H$.
\end{enumerate}
\renewcommand{\theenumi}{\arabic{enumi}}
\end{theorem}

Nguyen~\cite{nguyen} resolved the problem of characterising when a matroid has a non-trivial erection, and gave an explicit construction of the free erection (defined in Theorem~\ref{changethm}).
Roberts~\cite{roberts} exhibited a different procedure to construct all erections and from which the automorphism groups were characterised.


Independently, Knuth~\cite{knuthrandom} gave an algorithm for constructing all matroids on a finite set.
Within the algorithm, a certain procedure takes the copoints of a matroid $M$, a collection
of random sets and produces an erection of $M$ relative to these random sets.
If no random sets are added then the resulting erection is the free erection.

An attractive aspect of Knuth's construction is that the free erection is easily accessible (one simply adds no random sets).
The relation of the free erection to the conjectures above will be made clear in Theorem~\ref{lcequiv}.

We paraphrase Knuth's results in the next theorem.
For a collection $\T$ of subsets of $\sn$, let $\Expand(\T) := \{A\cup \{a\}: A\in\T \mbox{ and }a \in\sn-A\}$
and $\Random(\T):=\{A\subseteq S_n\;:\; A\supset T\mbox{ for all }T \in \T\}$,
essentially the strict filter generated by $\T$. 

Let $\Refine(\T,\T ')$, where $\T'$ is an arbitrary collection of sets,  be the collection resulting from: {Repeat until $A \cap B \subseteq C$ for some $C \in \T '$ for all $A,B \in \T$:
        if $A,B \in \T$, $A \neq B$ and $A \cap B \not\subseteq C$
        for any $C \in \T '$ then replace $A$ and $B$ in $\T$ by the set $A \cup B$.}

\begin{theorem}[Knuth\cite{knuthrandom}]\label{changethm} 
Let $M$ be a simple matroid on $S_n$ with rank function $\rho$ and collection of copoints $\H_M$.
\begin{enumerate}
\renewcommand{\theenumi}{\roman{enumi}}
\item If $N$ is an erection of $M$ then there exists a clutter $\A \subseteq \Random(\H_M)$, not necessarily unique and possibly empty,
such that
\begin{eqnarray*}
\F_{\rho(M)}(N) &=& \Refine(\Expand(\H_M)\cup \A, \H_M).
\end{eqnarray*}
\item Let $\F_{\rho(M)}(\Free(M)) = \Refine(\Expand(\H_M), \H_M)$. 
If \begin{eqnarray*}\F_{\rho(M)}(\Free(M)) &=& \{S_n\}\end{eqnarray*} then $\Free(M)=M$, otherwise the
copoints of the free erection are $\F_{\rho(M)}(\Free(M))$ and $\F_{\rho(M)+1}(\Free(M)):=\{S_n\}$.
\end{enumerate}
\renewcommand{\theenumi}{\arabic{enumi}}
\end{theorem}

It was also proven in ~\cite{knuthrandom} that the order of the replacements within the Refine function make no difference.
Formally this may be stated as:
\begin{eqnarray*}\Refine(\A,\T) &=& \Refine(\A - \{A,B\} \cup \Refine(\{A,B\},\T),\T)\end{eqnarray*}
and will be used implicitly in the next section.

Notice that if the free erection of a matroid $M$ is the trivial erection, then $M$ has no other erections and 
the ranks of $\Free(M)$ and $M$ are the same.
Otherwise, if the free erection of $M$ is not the trivial erection, then the free erection has rank $\rho(M)+1$, as do all other erections of
$M$ except the trivial erection.


\section{Inequalities concerning the free erection}

The following theorem shows that the free erection of a matroid $M$ 
is the matroid which contains the largest number of copoints amongst all erections of $M$.

\begin{theorem}
\label{theorem1}
Let $M$ be a simple rank-$r$ matroid and $N$ an erection of $M$, then $W_r(N) \leq W_r(\Free(M))$.
\end{theorem}

\begin{proof}
Let $\F_{r-1}=\F_{r-1}(M)$ be the copoints of $M$ 
and let $\R_{r-1} \subseteq\Random(\F_{r-1})$ be the collection of subsets of $S_n$ such that 
$\F_r(N) = \Refine(\Expand(\F_{r-1})\cup \R_{r-1},\F_{r-1})$.
Since $N$ is an erection of $M$, the existence of such a collection $\R_{r-1}$ is guaranteed due to Theorem~\ref{changethm}.

Let $\Expand(\F_{r-1})=\{E_1,\ldots ,E_m\}$ and define $E_i':=E_i$ for all $1\leq i \leq m$.
Notice that if $R_1,R_2 \in \R_{r-1}$ and $R_1,R_2 \supset E_i$, then 
\begin{eqnarray*}\Refine(\{E_i,R_1,R_2\},\F_{r-1}\} = \{R_1\cup R_2\}.\end{eqnarray*}
If there are more than two such $R$'s, 
then the $\Refine$ operation results in the single-element set containing their union.
This permits us to do the following:
for $i=1$ to  $m$, if $X \in \R_{r-1}$ and $X \supset E_i$ then let $E_i':= E_i'\cup X$ and remove $X$ from $\R_{r-1}$.
Repeat the previous operation until $\R_{r-1}$ is empty. 
Consequently, $E_i \subseteq E_i '$ for all $1\leq i \leq m$. 
For distinct $E_i,E_j \in \Expand(\F_{r-1})$, if 
$\Refine (\{E_i,E_j\},\F_{r-1}) = \{E_i \cup E_j\}$ then $\Refine (\{E_i' ,E_j'\},\F_{r-1}) = \{E_i' \cup E_j'\}$ because $E_i \cap E_j \subseteq E_i' \cap E_j'$.

Hence $\Refine(\{E_1,\ldots ,E_m\},\F_{r-1}\}= \Refine(\Expand(\F_{r-1}),\F_{r-1})$ has at least 
as many sets as 
\begin{eqnarray*}
\Refine(\{E_1',\ldots ,E_m'\},\F_{r-1}) &=& \Refine(\Expand(\F_{r-1})\cup \R_{r-1},\F_{r-1}))\\
&=&  \F_r(N).
\end{eqnarray*} 
\end{proof}

An alternative proof may be given using the results and terminology of Crapo~\cite{crapo70i},
in which it was shown that the lattice of all erections of a matroid $M$ contains a least element $\Free(M)$, 
and from which all other erections of $M$ may be obtained by partitioning (an anti-chain 
with respect to another anti-chain according to certain covering relations).

It would be misleading to think of the number of copoints in the free erection as possessing a monotone property.
The next example shows how slightly changing the copoints results in a completely different free erection.
Intuition would suggest the opposite to happen, however this is not the case.

\begin{example}
\label{example32}
Brylawski~\cite[p.171]{brylawski} gave the following example in which, 
by removing four copoints from a matroid and replacing them by their union, 
the number of copoints in the free erection actually increases.
Let $M_1$ and $M_2$ be the matroids on $S_8$ with 
\begin{eqnarray*}
\F_3(M_1) &=& \{ 123,124,1256,127,128,134,1357,136,138,1458,146,147,\\ 
	&& 167,168,178, 234,235,2367,238,245,246,247,248,257,258,\\ 
	&& 268,278,345,346,3478,356,358,368,456,457,467,468,567,\\ 
	&& 568,578,678\},\\
\F_3(M_2) &=& \F_3(M_1) - \{246,248,268,468\}\cup\{2468\}.
\end{eqnarray*}
One finds that $\F_4(\Free(M_1)) = \{S_8\}$ whereas
\begin{eqnarray*}
\F_4(\Free(M_2)) &=& \{ 123567, 1234, 1238, 1278, 124568, 1247, 134578, 1346,\\ 
	&& 1368, 1467, 1678, 2345, 234678, 2358, 2457, 2578, 3456,\\ 
	&& 3568, 4567, 5678 \}.
\end{eqnarray*}
\end{example}

In Dukes~\cite{dukes.slc} it was shown that the points-lines-planes conjecture is equivalent to bounding the number of copoints in the free erection of
a rank-3 matroid. 
We now generalize this to show the relation of free erection to the log-concavity conjecture for the Whitney numbers of
a matroid.
In this sense, the log-concavity conjecture may be considered a copoints-colines-coplanes conjecture.

\begin{theorem}
\label{lcequiv}
The log-concavity conjecture (\ref{lcstatement})(i) is true
if and only if for all simple matroids $M$,
\begin{eqnarray}
\label{wintersday}
W_{\rho(M)-1}(M)^2 \geq W_{\rho(M)-2}(M) W_{\rho(M)}(\Free(M)).
\end{eqnarray}
\end{theorem}

\begin{proof}
First assume the conjecture to be true and let $N=\Free(M)$ for some simple matroid $M$. Since $N\neq M$ (for otherwise it is trivial),
\begin{eqnarray*}
W_{\rho(N)-2} (N) ^2 & \geq & W_{\rho(N)-3}(N) W_{\rho(N)-1}(N)
\end{eqnarray*}
and $\F_{\rho(N)-3}(N)=\F_{\rho(M)-2}(M)$, $\F_{\rho(N)-2} (N) = \F_{\rho(M)-1}(M)$, $\F_{\rho(N)-1}(N) = \F_{\rho(M)}(\Free(M))$,
the inequality (\ref{wintersday}) holds.

Conversely, assume that for all simple matroids the inequality (\ref{wintersday}) is true. Suppose there
exists a matroid $N$ such that $W_i(N)^2 < W_{i-1}(N) W_{i+1}(N)$ for some $2<i<\rho(N)-1$.
Let $M=U_i(N)$ and from Theorem~\ref{theorem1}, since $\F_{i+1}(N)$ are the copoints of an erection of $M$,
\begin{eqnarray*}
W_{\rho(M)-1}(M)^2 &=& W_i(N)^2 ,\\
&<& W_{i-1}(N) W_{i+1}(N) ,\\
&=& W_{\rho(M)-2}(M) W_{\rho(M)}(\Free(M)),
\end{eqnarray*}
contradicting Equation (\ref{wintersday}).
\end{proof}

One might be tempted to insert the clause `such that $\Free(M)\neq M$' in the statement of the previous theorem, 
however this is not necessary due to Lemma~\ref{endclause}.

Nguyen's~\cite{nguyen} construction of the free erection can be viewed as wrapping 
the collection $\Expand(\F_{r-1})$ into a more compact collection 
$\Pair(\F_{r-1},\F_{r-2})$ to the effect that 
\begin{eqnarray*}
\Refine(\Expand(\F_{r-1}),\F_{r-1}) &=& \Refine(\Pair(\F_{r-1},\F_{r-2}),\F_{r-1})
\end{eqnarray*}
whereby $\Pair(\X,\Y):=\{X_1 \cup X_2 : X_1,X_2 \in \X \mbox{ and } X_1,X_2 \supset Y \in \Y\}$. 
So $\Pair(\F_{r-1},\F_{r-2})$ is the set containing the unions of pairs of copoints which have a common coline.
We give a short proof of the equivalence of the two collections under the $\Refine$ operation.

\begin{proposition}
\label{prop1}
Let $M$ be a simple rank-$r$ matroid on $S_n$ with copoints $\F_{r-1}$ and colines $\F_{r-2}$. Then
\begin{eqnarray*}
\Refine(\Expand(\F_{r-1}),\F_{r-1}) &=& \Refine(\Pair(\F_{r-1},\F_{r-2}),\F_{r-1}).
\end{eqnarray*}
\end{proposition}

\begin{proof}
The first point to note is that in the definition of the $\Refine$ function, 
the requirement that two sets be removed from $\Expand(\F_{r-1})$, if their union is to be included, is unnecessary. 
Should the two sets remain contained, then the $\Refine$ function will absorb both into their union as a final step.

Suppose $A,B \in \F_{r-1}$ and $A \cap B\in \F_{r-2}$. 
The existence of these pairs is guaranteed since for all $C \in \F_{r-2}$, 
the sets $\{ X-C: X \in \F_{r-1} \mbox{ and }X\supset C\}$ partition the set $S_n -C$.

Choose $a \in A\backslash B$ and $b \in B\backslash A$. 
Then $\{a\} \cup B, \, \{b\} \cup A \,\in\, \Expand(\F_{r-1})$
and $(\{a\} \cup B) \cap (\{b\} \cup A)\, =\, \{a,b\} \cup (A \cap B)$, 
but the set $\{a,b\} \cup (A \cap B)$ is not contained in any $C \in \F_{r-2}$ because $A\cap B \in\F_{r-2}$.
So the $\Refine$ function will replace the sets $\{a\} \cup B,\{b\}\cup A$ by the set $A \cup B$. 

Let us suppose this removal does not occur until all 
sets $A,B \in \F_{r-1}$ such that $A \cap B \in \F_{r-2}$ have been identified and their union inserted.
Thus $\Expand(\F_{r-1})$ contains all sets $A\cup B$ such 
that $A\cap B \in \F_{r-2}$ and $A\cap B \in \F_{r-2}$, 
which is precisely the collection $\Pair(\F_{r-1},\F_{r-2})$.
Now as the final operation, absorb all sets of the type 
$\{a\} \cup B \in \Expand(\F_{r-1})$ into the unions of $\Pair(\F_{r-1},\F_{r-2})$.
\end{proof}

\begin{example}
The Whitney numbers of a paving matroid are log-concave. 
A rank-$(r+1)$ paving matroid on $S_n$ is the free erection of $N_r$, the $r$-truncation 
of the boolean algebra on $S_n$ (i.e. $N_r$ is the matroid with copoints ${S_n \choose r-1}$.)
Thus 
\begin{eqnarray*}
W_i(M) && \left\{ 
	\begin{array}{ll}
	= \; W_i(N) \;= \;{n \choose i}, & \mbox{ if } i \leq r-1, \\
	\leq \; W_{i-1} (\Free(N_r)) \;= \;{n\choose i}, &  \mbox{ if } i= r, \\
	= \; 1, & \mbox{ if } i \;= \;r+1.
	\end{array}
	  \right.
\end{eqnarray*}
The outstanding case is $W_{r}(M)^2 \geq W_{r-1}(M) W_{r+1}(M)$, but this
follows immediately from Lemma~\ref{endclause}.
Hence the sequence of Whitney numbers is log-concave. 
\end{example}

For a simple matroid $M$ on $S_n$ and $\pi$ a permutation of $S_n$,
let $M(\pi)$ be the matroid obtained from $M$ by permuting all elements of $M$ by $\pi$ in the
natural way.
We present the following bound on the number of copoints in the free erection.
The bound is very case specific, as will be seen in the examples mentioned after the proof.
For matroids whose free erection is the trivial erection the bound becomes
an equality, unsurprisingly, if there exist two copoints containing a common coline,
whose union is the ground set of the matroid. 

\begin{theorem}
\label{freebound}
Let $M$ be a simple rank-$r$ matroid on $S_n$. For each coline $F$, let
$a_F(M) := |\{X \in \F_{r-1}(M): \min_M(X) = F\}|$ where $\min_M(X)$ denotes
the lexicographically smallest coline contained in the copoint $X$. Then
\begin{eqnarray}
\label{philips}
W_r(\Free(M)) & \leq & \min_{\pi} \sum_{F \in \F_{r-2}(M(\pi))} {a_F(M(\pi)) \choose 2}.
\end{eqnarray}
\end{theorem}

\begin{proof}
\newcommand{\mypair}{\mbox{\scshape{Pair}}^{\star}}
From Proposition~\ref{prop1},
\begin{eqnarray*}
\F_r(\Free(M)) &=& \Refine(\Pair(\F_{r-1},\F_{r-2}),\F_{r-1}).
\end{eqnarray*}
Let $F(X_1,X_2)$ be the unique copoint containing the two distinct colines $X_1,X_2$ of $M$ and
let $\min_M(X)$ denote the lexicographically smallest coline containing the copoint $X$. 
Define 
\begin{eqnarray*}
\mypair(M) &:=& \{\mbox{distinct copoints }(X,Y)\mbox{ of } M\;:\; \min_M(X)=\min_M(Y) = F,\\ && \hspace*{14em} \; \mbox{for all colines }F \mbox{ of }M\}.
\end{eqnarray*}

Given $H\in \F_r(\Free(M))$, let $\alpha_1(H) < \alpha_2(H) < \alpha_3(H)$ be the three lexicographically smallest 
colines of $M$ contained in $H$. 
For every such $H\in \F_r(\Free(M))$, the pair of copoints $(F(\alpha_1(H),\alpha_2(H)),F(\alpha_1(H),\alpha_3(H)))$ (of $M$) uniquely determine $H$.
Thus $|\F_r(\Free(M))| \leq |\mypair(\F_{r-1},\F_{r-2})|$.
To evaluate $\mypair$ notice that if
\begin{eqnarray*}
a_F(M) &=& |\{\mbox{copoints } (X,Y)\mbox { of }M: \min_M(X) = \min_M(Y) = F\}|,
\end{eqnarray*}
then
\begin{eqnarray*}
|\mypair(\F_{r-1},\F_{r-2})| &=& \sum_{F\in\F_{r-2}(M)} {a_F(M) \choose 2}.
\end{eqnarray*}
(Of course $\sum_{F\in\F_{r-2}} a_F(M) = |\F_{r-1}(M)|$.)
Since the number of copoints in the free erection is invariant under permutations of the ground set, 
we can take the minimum over all such permutations.
\end{proof}

\begin{example} The varying performance of the bound may be seen through the following examples;
\label{example36}
\renewcommand{\theenumi}{\roman{enumi}}
\begin{enumerate}
\item
Let $M_1$ and $M_2$ be the matroids in Example~\ref{example32}.
We have that $W_4(\Free(M_1))=1$ and $W_4(\Free(M_2)) = 20$.
The permutation $\pi=12456387$ (not unique) minimizes both sums in (\ref{philips}), 
which are 32 and 30 for $M_1$ and $M_2$, respectively.
\item Let $M_3$ be the rank-4 matroid on $S_8$ with $\F_3(M_3) = {S_8 \choose 3}$.
Clearly $M(\pi) = M$ for all permutations $\pi$ so the minimum in (\ref{philips}) is $\sum_{1}^5 i{7-i\choose 2} = 70$. 
The new copoints are $\F_4(\Free(M_3)) = {S_8 \choose 4}$, giving $W_4(\Free(M_3)) = 70$, and there is equality in (\ref{philips}).
\item
Let $M_4$ be the rank-3 matroid on $S_8$ with
\begin{eqnarray*}\F_2(M_4) &=& \{12,13,14,15,16,16,17,18,2345678\}.\end{eqnarray*}
The minimum is attained at any permutation $\pi$ with $\pi_1=8$, so that  
$\F_2(M_4(\pi)) = \{1234567,18,28,38,48,58,68,78\}$
and $a_{\{1\}}(M_4(\pi)) = 2$, $a_{\{2\}}(M_4(\pi)) = \ldots = a_{\{8\}}(M_4(\pi)) = 1$ and so the sum in (\ref{philips}) is 1.
Notice $\Free(M_4)= M_4$ and again there is equality in (\ref{philips}).
\item Let $M_5$ be the rank-3 matroid on $S_8$ with copoints (from Nguyen~\cite[Example 1]{nguyen})
\begin{eqnarray*}
\F_2(M_5) &=& \{123,14,15,16,17,18,24,258,26,27,34,35,\\ && 368,37,45,46,478,56,57,67\}.
\end{eqnarray*}
The copoints of the free erection are 
\begin{eqnarray*}
\F_4(\Free(M_5)) &=& \{1234,123568,1237,145,146,1478,157,\\ &&167,24578,246,267,345,34678,357,\\ && 456,567\}
\end{eqnarray*}
Thus $W_4(\Free(M_5)=16$ and the minimum of the sum is 26, attained at $\pi=18234567$.
\item 
Let $M_6$ be the rank-3 matroid on $S_9$ with copoints (from Nguyen~\cite[Example 2]{nguyen})
\begin{eqnarray*}
\F_2(M_6) &=& \{123,14,156,17,18,19,24,25,26,279,\\ && 28,34,35,36,37,38,39,45,46,47,48,49,\\ &&57,58,59,678,69,89\}.
\end{eqnarray*}
The free erection of $M_6$ is the trivial erection so $W_4(\Free(M_6))=1$ and the minimum of the sum is 44,
attained at $\pi=12678459$.
\end{enumerate}
\renewcommand{\theenumi}{\arabic{enumi}}
\end{example}

\section{The free erection of a rank-3 geometry}
The definition of the free erection, both in the present paper and previous papers 
addressing the subject ~\cite{nguyen,revjakin,roberts,crapo70i},
is the outcome of a stringent sequence of joining operations on the copoints.
Attempts to describe the free erection any further lead to either a reformulation of the $\Refine$ function, 
or a new collection $\A$ formed from copoints $\F_r$ such that
the free erection is $\Refine(\A,\F_r)$.
Indeed, in what is described below we do not achieve any further insights into the free erection.
However, the formulation below has the benefit that, in the rank-3 case, one can `visualize' the free erection,
thereby making it accessible to graph theorists. 

The first class of matroids with non-trivial erections are rank-3 matroids.

The relationship of free erections of rank-3 matroids to the points-lines-planes conjecture was shown in \cite{dukes.slc} and we restate it here.
The points-lines-planes conjecture states that for all simple rank-4 matroids $M$ on $S_n$,
\begin{eqnarray*}
W_2(M)^2 &\geq& \dfrac{3(n-1)}{2(n-2)} W_1(M)W_{3}(M).
\end{eqnarray*}
Recall that the points-lines-planes conjecture is a statement about matroids of rank (at least) 4, whereas
the statement of the next theorem concerns matroids of rank-3.
\begin{theorem}[Dukes{~\cite[Theorem 4.4]{dukes.slc}}]
\label{myequiv}
{The points-lines-planes conjecture is true if and only if it is true for the
class of rank-4 matroids that are the free erection of some rank-3 matroid.}
\end{theorem}


The copoints of a simple rank-3 matroid are a 2-partition of the set $S_n$ (see for example Aigner~\cite[p.258]{aigner}).
We use the notation of Diestel~\cite{diestel}.
Let us say a bipartite graph 
$G(V,E)$ with bipartition $\{V_1,V_2\}$
has property $\property$ if 
\begin{itemize}
\item $V_1$ and $V_2$ are non-empty,
\item for all $v_1,v_1' \in V_1$ there exists a unique $v_2\in V_2$ such that\\ $(v_1,v_2),(v_1',v_2)\in E(G)$,
\item $\delta(G)\geq 2$.
\end{itemize}

Bipartite graphs with property $\property$ are in one-to-one correspondence with simple rank-3 matroids (up to isomorphism) since
they represent the covering relations of the colines to copoints.
In the rank 3 case the colines are the points and the copoints are the lines.
For $v \in V_2$, let $\preimage{v}=\{v_1 \in V_1: (v_1,v)\in E(G)\}$ and similarly for $A \subseteq V_2$ let 
$\preimage{A}$ be the union of the sets $\preimage{v}$, for all $v \in A$.
Denote by $\myrest{G}{A\cup\preimage{A}}$ the subgraph of $G$ induced by restriction to the vertex set $A\cup \preimage{A}$.

Define $\property(G):=|\{ A\subset V_2: \myrest{G}{A\cup \preimage{A}} \mbox{ has property } \property\}|$.
Since the number of copoints in the free erection (for matroids which have a non-trivial erection) is precisely the number $\property(G)$, we have

\begin{proposition}
\label{plpprop}
The points-lines-planes conjecture is true if and only if for all bipartite graphs 
$G(V,E)$ 
with property $\property$
\begin{eqnarray*}
\property(G) & \leq & \dfrac{2m^2(n-2)}{3n(n-1)},
\end{eqnarray*}
where $n=|V_1|$, $m=|V_2|$ and $\{V_1,V_2\}$ is the bipartition of $G$.
\end{proposition}

\begin{example}
Consider the rank-3 matroid $M_5$ as mentioned in Example~\ref{example36}.
Figure 1 shows the covering relations between the copoints and colines of $M_5$, thereby
defining the associated bipartite graph $G_M$.
The vertices of $G_M$ are $V_1=\{1,2,\ldots ,8\}$ and $V_2=\{123,14,\ldots,67\}$.

If we choose $A=\{24,258,27,45,478,57\}$ then we find that the subgraph $\myrest{G}{A\cup\preimage{A}}$ induced by 
the vertices $\{24,258,27,45,478,57,2,4,5,7,8\}$ has property $\property$,
hence $\preimage{A}=\{24578\}\in \F_4(M_5)$ (see Figure 2).

However, if $A=\{258,26,368,56\}$ then we find that the subgraph $\myrest{G}{A\cup\preimage{A}}$ induced by
the vertices $\{258,26,368,56,2,3,5,6,8\}$ is such that there does not exist $v$ with both $(2,v),(3,v) \in E(\myrest{G}{A\cup\preimage{A}})$,
and similarly for the pair $\{3,5\} \in V_1$. So $\myrest{G}{A\cup\preimage{A}}$ does not have property $\property$ and
hence $\preimage{A}=\{23568\}\not\in \F_4(M_5)$ (see Figure 3).

\begin{figure}
\scalebox{0.7}{\includegraphics{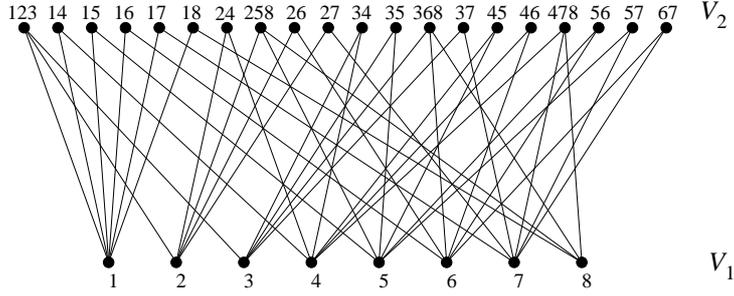}}
\caption{The bipartite graph 
$G(\F_1(M_5) \cup \F_2(M_5),E)$
representing the colines (points) and copoints (lines) of $M_5$ in Example~\ref{example36}(v)}
\end{figure}
\begin{figure}
\scalebox{0.7}{\includegraphics{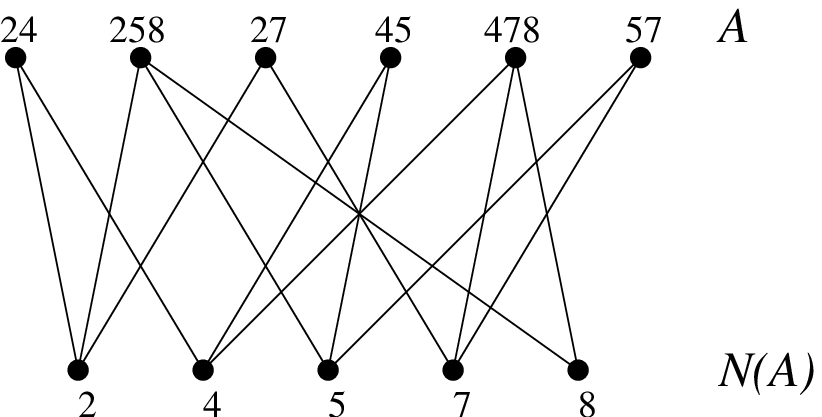}}
\caption{The induced subgraph of the bipartite graph of $M_5$ when restricted to $A\cup\preimage{A} = \{24,258,27,45,478,$ $57,2,4,5,7,8\}$}
\end{figure}
\begin{figure}
\scalebox{0.7}{\includegraphics{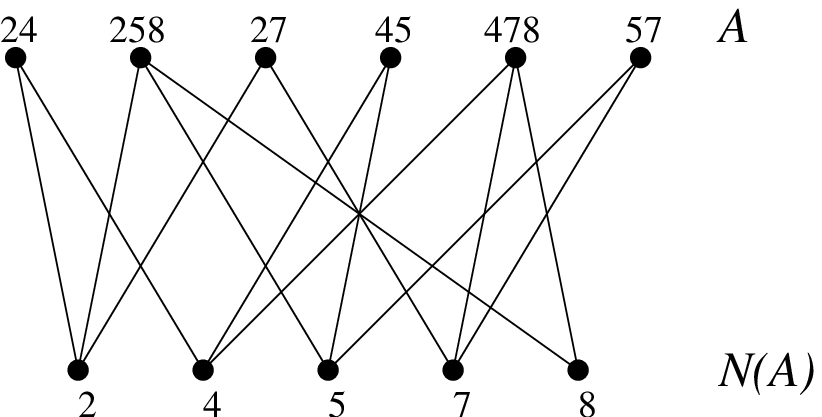}}
\caption{The induced subgraph of the bipartite graph of $M_5$ when restricted to $A\cup\preimage{A} = \{258,26,368,56,$ $2,3,5,6,8\}$. }
\end{figure}
\end{example}

A tempting case to first resolve is the situation of Kung~\cite{kung}, in which every line contains an equal number of points ($t$ say).
In the setting of Proposition~\ref{plpprop} we find that $d(v)=(n-1)/(t-1)$ for $v \in V_1$ and $d(v)=t$ for $v \in V_2$. 
This is the Steiner system $S(2,t,n)$ with $m=n(n-1)/t(t-1)$ lines (blocks).
As mentioned in the introduction, the cases for $t=3,4$ have been resolved.
It is possible that the theory of Steiner system can be applied to this case in Proposition~\ref{plpprop} to resolve the conjecture for $t\geq 5$.


\begin{thebibliography}{99}
\bibitem{aigner} M. Aigner, Combinatorial Theory, Springer-Verlag, Berlin/Heidelberg/New York, 1979.
\bibitem{brylawski} T. Brylawski, `Constructions' in {\it{Theory of Matroids}} (Neil White, ed.), Cambridge University Press, 1986.
\bibitem{crapo70i} H. H. Crapo, Erecting Geometries, {\it{Proc. 2$^{nd}$ Chapel Hill Conf. on Comb. Math.}} (1970), 74--99.
\bibitem{crapo70ii} H. H. Crapo, Erecting geometries, {\it{Ann. New York Acad. Sci.}} {\bf{175}} (1970), 89--92.
\bibitem{diestel} R. Diestel, Graph Theory, Springer-Verlag, New York, 1997.
\bibitem{dukesaus} W. M. B. Dukes, Bounds on the number of generalized partitions and some applications, {\it{Australas. J. Comb.}} {\bf{28}} (2003), 257--262.
\bibitem{dukes.slc} W. M. B. Dukes, On the number of matroids on a finite set, {\it{S\'{e}m. Lothar. Combin.}} (electronic) {\bf{51}} (2004), to appear.
\bibitem{knuthbound} D. E. Knuth, The Asymptotic Number of Geometries, {\it{J. Combin. Theory A}} {\bf{16}} (1974), 398--400.
\bibitem{knuthrandom} D. E. Knuth, Random Matroids, {\it{Discrete Math.}} {\bf{12}} (1975), 341--358.
\bibitem{kung} J. P. S. Kung, On the lines-planes inequality for matroids, {\it{J. Combin. Theory Ser. A}} {\bf{91}} (2000), 363--368.
\bibitem{lasvergnas} M. Las Vergnas, On certain constructions for matroids, {\it{Proc. 5th Comb. Conf., Congressus Numerantium, No. XV, Utilitas Math., Win. Man.}} (1976) 395--404.
\bibitem{mason} J. H. Mason, Matroids: Unimodal conjectures and Motzkin's theorem, {\it{Combinatorics}} (D. J. A. Welsh and D. R. Woodall, eds.), Institute of Math. and Appl. (1972), 207--221.
\bibitem{nguyen} H. Q. Nguyen, Constructing the free erection of a geometry, {\it{J. Combin. Theory B}} {\bf{27}} (1979), 216--224.
\bibitem{revjakin} A. M. Revjakin, Erections of combinatorial geometries, {\it{Vestnik. Moskov. Univ. Ser. I Mat. Meh.}} {\bf{31}}(4) (1976), 59--62.
\bibitem{roberts} L. Roberts, All erections of a combinatorial geometry and their automorphism groups, {\it{Lecture Notes in Mathematics}} {\bf{452}} (1975), 210--213.
\bibitem{rota64} G. C. Rota, Combinatorial Theory, old and new, {\it{Actes, Congres. Intern. Math.}} (1970) Tome 3, 229--233.
\bibitem{seymour} P. D. Seymour, On the points-lines-planes conjecture, {\it{J. Combin. Theory B}} {\bf{33}} (1982), 17--26.
\bibitem{stone} J. R. Stonesifer, Logarithmic concavity for edge lattices of graphs, {\it{J. Combin. Theory A}} {\bf{18}} (1975), 36--46.
\bibitem{welshbook} D. J. A. Welsh, Matroid Theory, Academic Press, 1976.
\end{thebibliography}
\end{document}